\documentclass[12pt]{amsart}
\usepackage{epsfig, graphics, psfrag, color}
\usepackage{amscd} 
\newcommand{\D}{{\partial}} % Boundary
\newcommand{\HH}{{\mathbb{H}}}
\newcommand{\RR}{{\mathbb{R}}}
\newcommand{\ZZ}{{\mathbb{Z}}}

\newcommand{\CC}{{\mathbb{C}}}

\def\co{\colon\thinspace}

\newcommand{\F}{{\mathcal{F}}}
\newcommand{\R}{{\mathcal{R}}}

\theoremstyle{plain}
\newtheorem{theorem}{Theorem}[section]

\newtheorem{cor}[theorem]{Corollary}
\newtheorem{lemma}[theorem]{Lemma}
\newtheorem{prop}[theorem]{Proposition}

\newtheorem*{namedtheorem}{\theoremname}
\newcommand{\theoremname}{testing}
\newenvironment{named}[1]{\renewcommand{\theoremname}{#1}\begin{namedtheorem}}{\end{namedtheorem}}

\theoremstyle{definition}
\newtheorem{define}[theorem]{Definition}

\newtheorem{example}[theorem]{Example}
\newtheorem*{remark}{Remark}

\begin{document}

%\hsize=6truein
%\hoffset=.3truein

\title{The length of unknotting tunnels}

\author{Daryl Cooper}
\author{Marc Lackenby}
\author{Jessica S. Purcell}

%% \address{Department of Mathematics, University of California, Santa
%% Barbara, Santa Barbara, CA 93106, USA}

%% \address{Mathematical Institute, 24--29 St Giles', Oxford, OX1 3LB,
%% 	England} 

%% \address{Department of Mathematics, Brigham Young University, Provo,
%% 	UT 84604, USA}

%\thanks{ \today}

\begin{abstract}
	We show there exist tunnel number one hyperbolic 3--manifolds with
	arbitrarily long unknotting tunnel.  This provides a negative answer
	to an old question of Colin Adams.
\end{abstract}

\maketitle

\newcommand{\mat}[2][cccc]{\left(\begin{array}{#1} #2\\
	\end{array}\right)}

\section{Introduction}\label{sec:intro}
In a paper published in 1995 \cite{adams:tunnels}, Colin Adams studied
geometric properties of hyperbolic tunnel number one manifolds.  A
\emph{tunnel number one manifold} is defined to be a compact
orientable 3--manifold $M$ with torus boundary component(s), which
contains a properly embedded arc $\tau$, the exterior of which is a
handlebody.  The arc $\tau$ is defined to be an \emph{unknotting
tunnel} of $M$.

When a tunnel number one manifold $M$ admits a hyperbolic structure,
there is a unique geodesic arc in the homotopy class of $\tau$.  If
$\tau$ runs between distinct boundary components, Adams showed that
its geodesic representative has bounded length, when measured in the
complement of a maximal horoball neighborhood of the cusps.  He asked
a question about the more general picture: does an unknotting tunnel
in a hyperbolic 3--manifold always have bounded length?

In response, Adams and Reid showed that when the tunnel number one
manifold is a 2--bridge knot complement, that unknotting tunnels have
bounded length \cite{adams-reid}.  Akiyoshi, Nakagawa, and Sakuma
showed that unknotting tunnels in punctured torus bundles actually
have length zero \cite{ans}, hence bounded length.

Sakuma and Weeks also studied unknotting tunnels in 2--bridge knots
\cite{sakuma-weeks}.  They found that any unknotting tunnel of a
2--bridge knot was isotopic to an edge of the canonical polyhedral
decomposition of that knot, first explored by Epstein and Penner
\cite{epstein-penner}.  They conjectured that all unknotting tunnels
were isotopic to edges of the canonical decomposition.  Heath and Song
later showed by example that not all unknotting tunnels could be
isotopic to edges of the canonical decomposition \cite{heath-song}.
However, the question of whether unknotting tunnels have bounded
length remained unanswered.

In this paper we finally settle the answer to this question.  We show
that, in fact, the answer is no.  There exist tunnel number one
manifolds with arbitrarily long unknotting tunnel.

\begin{named}{Theorem \ref{thm:long}}
	There exist finite volume one--cusped hyperbolic tunnel number one
	manifolds for which the geodesic representative of the unknotting
	tunnel is arbitrarily long, as measured between the maximal horoball
	neighborhood of the cusp.
\end{named}

Note we are not claiming here that the unknotting tunnel in these
examples is ambient isotopic to a geodesic.  Such examples can in fact
be constructed, but the argument is more complex and will appear in a
companion paper \cite{lackenby-purcell-2}.  However, note that Theorem
\ref{thm:long} does force the unknotting tunnels in these examples to
be arbitrarily long, because the length of a properly embedded arc is
at least that of the geodesic in its homotopy class.

We prove Theorem \ref{thm:long} in two ways.  The first proof, which
appears in Section \ref{sec:long-tunnels}, is geometric and partially
non-constructive.  We analyze the infinite--volume hyperbolic
structures on the compression body $C$ with negative boundary a torus,
and positive boundary a genus 2 surface.  A guiding principle is that
geometric properties of hyperbolic structures on $C$ should often have
their counterparts in finite--volume hyperbolic 3--manifolds with
tunnel number one.  For example, any geometrically infinite hyperbolic
structure on $C$ is the geometric and algebraic limit of a sequence of
geometrically finite hyperbolic structures on $C$, and it is also the
geometric limit of a sequence of finite--volume hyperbolic
3--manifolds with tunnel number one.  It is by finding suitable
sequences of hyperbolic structures on $C$ that Theorem \ref{thm:long}
is proved.  In particular, the proof gives very little indication of
what the finite--volume hyperbolic 3--manifolds actually are.

The geometric proof of Theorem \ref{thm:long} leads naturally to the
study of geometrically finite structures on the compression body $C$
and their geometric properties.  We include some background material
in Sections \ref{sec:prelim} and \ref{sec:ford}.  However, we postpone
a more extensive investigation of geometrically finite structures on
$C$ to a companion paper \cite{lackenby-purcell-2}.

The second proof is more topological, and appears in Section
\ref{sec:dehn}.  The idea is to start with a tunnel number one manifold
with two cusps.  An argument using homology implies that there exist
Dehn fillings on one cusp which yield a tunnel number one manifold whose
core tunnel must be arbitrarily long.

A consequence of the second proof is that the resulting tunnel number
one manifold cannot be the exterior of a knot in a homology sphere.  In
Section \ref{sec:homology}, we modify the construction of the first
proof to show there do exist tunnel number one manifolds with long
tunnel which are the exterior of a knot in a homology sphere.  It
seems likely that the Dehn filling construction in Section
\ref{sec:dehn} can be modified to produce hyperbolic knots in homology
spheres with long unknotting tunnels.  However, to establish this, a
substantially different method of proof would be required.

Although we construct examples of knots in homology 3--spheres with
long unknotting tunnels, we do not obtain knots in the 3--sphere using
our methods.  It would be interesting to determine whether such
sequences of knots exist.  If they do, can explicit diagrams of such
knots be found?

\subsection{Acknowledgements}
Lackenby and Purcell were supported by the Leverhulme trust.  Lackenby
was supported by an EPSRC Advanced Research Fellowship.  Purcell was
supported by NSF grant DMS-0704359.  Cooper was partially supported by
NSF grant DMS-0706887.

\section{Background and preliminary material}\label{sec:prelim}
In this section we will review terminology and results used throughout
the paper.

The first step in the proof of Theorem \ref{thm:long} is to show there
exist geometrically finite structures on a compression body $C$ with
arbitrarily long tunnel.  We begin by defining these terms.

\subsection{Compression bodies}

%% define compression body
A \emph{compression body} $C$ is either a handlebody, or the result of
taking a closed, orientable (possibly disconnected) surface $S$ cross
an interval $[0,1]$, and attaching 1--handles to $S\times\{1\}$.  The
\emph{negative boundary}, denoted $\partial_-C$, is $S \times \{ 0\}$.
When $C$ is a handlebody, $\partial_{-}C = \emptyset$.  The
\emph{positive boundary} is $\partial C \setminus \partial_{-}C$, and
is denoted $\partial_{+}C$.

%% define unknotting tunnel for compression body
Throughout this paper, we will be interested in compression bodies $C$
for which $\partial_{-}C$ is a torus and $\partial_{+}C$ is a genus
$2$ surface.  We will refer to such a manifold as a
\emph{$(1,2)$--compression body}, where the numbers $(1,2)$ refer to
the genus of the boundary components.

\begin{figure}
	\centerline{\includegraphics[width=3in]{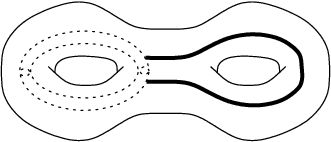}}
\caption{The $(1,2)$--compression body. The core tunnel is the thick
	line shown, with endpoints on the torus boundary.}
\label{fig:comp-body}
\end{figure}

Let $\tau$ be the union of the core of the attached 1--handle with
two vertical arcs in $S \times [0,1]$ attached to its endpoints.  Thus,
$\tau$ is a properly embedded arc in $C$, and $C$ is a regular
neighborhood of $\partial_- C \cup \tau$.  We refer to $\tau$ as the
\emph{core tunnel}.  See Figure \ref{fig:comp-body}.

%% review fundamental group of compression body, generators \alpha,
%% \beta, \gamma.
Note that the fundamental group of a $(1,2)$--compression body $C$ is
isomorphic to $(\ZZ\times\ZZ)\ast \ZZ$.  We will denote the generators
of the $\ZZ \times \ZZ$ factor by $\alpha$, $\beta$, and we will
denote the generator of the second factor by $\gamma$.

\subsection{Hyperbolic structures}

Let $C$ be a $(1,2)$--compression body.  We are interested in complete
hyperbolic structures on the interior of $C$.  We obtain a hyperbolic
structure on $C \setminus \D C$ by taking a discrete, faithful
representation $\rho\co \pi_1(C) \to {\rm PSL}(2,\CC)$ and considering the
manifold $\HH^3/\rho(\pi_1(C))$.

\begin{define}
A discrete subgroup $\Gamma < {\rm PSL}(2,\CC)$ is \emph{geometrically
finite} if $\HH^3/\Gamma$ admits a finite--sided, convex fundamental
domain.  In this case, we will also say that the manifold
$\HH^3/\Gamma$ is \emph{geometrically finite}.
\label{def:gf-group}
\end{define}

Geometrically finite groups are well understood.  In this paper, we
will often use the following theorem of Bowditch (and its corollary,
Corollary \ref{cor:bowditch} below).

\begin{theorem}[Bowditch, Proposition 5.7 \cite{bowditch}]
If a subgroup $\Gamma < {\rm PSL}(2,\CC)$ is geometrically finite, then
every convex fundamental domain for $\HH^3/\Gamma$ has finitely many
faces.
\label{thm:bowditch}
\end{theorem}

\begin{define}
For $C$ a $(1,2)$--compression body, we will say that a discrete, faithful
representation $\rho$ is \emph{minimally parabolic} if for all $g \in
\pi_1(C)$, $\rho(g)$ is parabolic if and only if $g$ is conjugate to
an element of the fundamental group of the torus boundary component
$\D_-C$.
\label{def:min-parabolic}
\end{define}

\begin{define}
A discrete, faithful representation $\rho\co\pi_1(C)\to {\rm PSL}(2,\CC)$ is
a \emph{minimally parabolic geometrically finite uniformization of
$C$} if $\rho$ is minimally parabolic, $\rho(\pi_1(C))$ is
geometrically finite as a subgroup of ${\rm PSL}(2,\CC)$, and
$\HH^3/\rho(\pi_1(C))$ is homeomorphic to the interior of $C$.
\label{def:gf}
\end{define}

%% PERHAPS SOME CONDITIONS IN ABOVE DEF ARE REDUNDANT?

It is a classical result, due to Bers, Kra, and Maskit (see
\cite{bers74}), that the space of conjugacy classes of minimally
parabolic geometrically finite uniformizations of $C$ may be
identified with the Teichm\"uller space of the genus 2 boundary
component $\D_+C$, quotiented out by ${\rm Mod}_0(C)$, the group of
isotopy classes of homeomorphisms of $C$ which are homotopic to the
identity.

In particular, note that the space of minimally parabolic
geometrically finite uniformizations is path connected.

\subsection{Isometric spheres and {F}ord domains}

The tool we use to study geometrically finite representations is that
of Ford domains.  We define the necessary terminology in this section.

Throughout this subsection, let $M=\HH^3/\Gamma$ be a hyperbolic
manifold with a single rank two cusp, for example, the
$(1,2)$--compression body.  In the upper half space model for $\HH^3$,
assume the point at infinity in $\HH^3$ projects to the cusp.  Let $H$
be any horosphere about infinity.  Let $\Gamma_\infty < \Gamma$ denote
the subgroup that fixes $H$.  By assumption, $\Gamma_\infty =
\ZZ\times\ZZ$.

%% define isometric sphere
\begin{define}
For any $g \in \Gamma \setminus\Gamma_\infty$, $g^{-1}(H)$ will be a
horosphere centered at a point of $\CC$, where we view the boundary at
infinity of $\HH^3$ to be $\CC \cup \{\infty\}$.  Define the set $S_g$
to be the set of points in $\HH^3$ equidistant from $H$ and
$g^{-1}(H)$.  $S_g$ is the \emph{isometric sphere} of $g$.
\label{def:isometric-sphere}
\end{define}

Note that $S_g$ is well--defined even if $H$ and $g^{-1}(H)$ overlap.
It will be a Euclidean hemisphere orthogonal to the boundary $\CC$ of
$\HH^3$.

At first glance, it may seem more natural to consider points
equidistant from $H$ and $g(H)$, rather than $g^{-1}(H)$ as in
Definition \ref{def:isometric-sphere}.  However, we use the historical
definition of isometric spheres in order to make use of the following
classical result, which we include as a lemma.  A proof can be found,
for example, in Maskit's book \cite[Chapter IV, Section G]{maskit}.

\begin{lemma}
For any $g\in\Gamma\setminus\Gamma_\infty$, the action of $g$ on
$\HH^3$ is given by inversion in $S_g$ followed by a Euclidean
isometry. \qed
\label{lemma:isosphere-invert}
\end{lemma}

The following is well known, and follows from standard calculations in
hyperbolic geometry.  We give a proof in \cite{lackenby-purcell-2}.

\begin{lemma}
If $$g=\mat{a&b\\c&d}\in {\rm PSL}(2,\CC),$$ then the center of the
Euclidean hemisphere $S_{g^{-1}}$ is $g(\infty) = a/c$.  Its Euclidean
radius is $1/|c|$.
\label{lemma:iso-center-rad}
\end{lemma}

%% define Ford domain

Let $B_g$ denote the \emph{open} half ball bounded by $S_g$.  Define
$\F$ to be the set
$$\F = \HH^3 \setminus \bigcup_{g\in \Gamma \setminus \Gamma_\infty}
B_g.$$ Note $\F$ is invariant under $\Gamma_\infty$, which acts by
Euclidean translations on $\HH^3$.

When $H$ bounds a horoball $H_\infty$ that projects to an embedded
horoball neighborhood about the rank 2 cusp of $M$, $\F$ is the set of
points in $\HH^3$ which are at least as close to $H_\infty$ as to any
of its translates under $\Gamma \setminus \Gamma_\infty$.  Such an
embedded horoball neighborhood of the cusp always exists, by the
Margulis lemma.

\begin{define}
A \emph{vertical fundamental domain for $\Gamma_\infty$} is a
fundamental domain for the action of $\Gamma_\infty$ cut out by
finitely many vertical geodesic planes in $\HH^3$.
\end{define}

\begin{define}
A \emph{Ford domain} of $M$ is the intersection of $\F$ with a
vertical fundamental domain for the action of $\Gamma_\infty$.
\label{def:ford-domain}
\end{define}

A Ford domain is not canonical because the choice of fundamental
domain for $\Gamma_\infty$ is not canonical.  However, for the
purposes of this paper, the region $\F$ in $\HH^3$ is often more
useful than the actual Ford domain.

%% Result (Bowditch):  In 3D, geometrically finite if and only if every
%% convex polyhedron has a finite number of faces.  Thus geometrically
%% finite if and only if Ford domain has a finite number of sides.
Note that Ford domains are convex fundamental domains.  Thus we have
the following corollary of Bowditch's Theorem \ref{thm:bowditch}.

\begin{cor}
$M=\HH^3/\Gamma$ is geometrically finite if and only if a Ford
domain for $M$ has a finite number of faces.
\label{cor:bowditch}
\end{cor}

\subsection{Visible faces and {F}ord domains}

\begin{define}
Let $g \in \Gamma\setminus\Gamma_\infty$.  The isometric sphere $S_g$
is called \emph{visible from infinity}, or simply \emph{visible}, if
it is not contained in $\bigcup_{h \in \Gamma\setminus(\Gamma_\infty
\cup \Gamma_\infty g)} \bar{B_h}$.  Otherwise, $S_g$ is called
\emph{invisible}.

Similarly, suppose $g, h \in \Gamma\setminus\Gamma_\infty$, and $S_g
\cap S_h \cap \HH^3$ is nonempty.  Then the edge of intersection $S_g
\cap S_h$ is called \emph{visible} if $S_g$ and $S_h$ are visible and
their intersection is not contained in $\bigcup_{k \in
\Gamma\setminus(\Gamma_\infty \cup \Gamma_\infty g \cup \Gamma_\infty
h)} \bar{B_k}$.  Otherwise, it is \emph{invisible}.
\label{def:visible}
\end{define}

The faces of $\F$ are exactly those that are visible from infinity.

In the case where $H$ bounds a horoball that projects to an embedded
horoball neighborhood of the rank 2 cusp of $M$, there is an
alternative interpretation of visibility.  An isometric sphere $S_g$
is visible if and only if there exists a point $x$ in $S_g$ such that
for all $h \in \Gamma \setminus (\Gamma_\infty \cup \Gamma_\infty g)$,
the hyperbolic distance $d(x, h^{-1}(H))$ is greater than the
hyperbolic distance $d(x, H)$.  Similarly, an edge $S_g \cap S_h$ is
visible if and only if there exists a point $x$ in $S_g \cap S_h$ such
that for all $k \in \Gamma \setminus (\Gamma_\infty \cup \Gamma_\infty
g \cup \Gamma_\infty h)$, the hyperbolic distance $d(x, H)$ is
strictly less than the hyperbolic distance $d(x, k^{-1}(H))$.

We present a result that allows us to identify minimally parabolic
geometrically finite uniformizations.

\begin{lemma}
Suppose $\rho\co \pi_1(C) \to {\rm PSL}(2,\CC)$ is a geometrically
finite uniformization.  Suppose none of the visible isometric spheres
of the Ford domain of $\HH^3/\rho(\pi_1(C))$ are visibly tangent on
their boundaries.  Then $\rho$ is minimally parabolic.
\label{lemma:min-parabolic}
\end{lemma}

By \emph{visibly} tangent, we mean the following.  Set $\Gamma =
\rho(\pi_1(C))$, and assume a neighborhood of infinity in $\HH^3$
projects to the rank two cusp of $\HH^3/\Gamma$, with $\Gamma_\infty <
\Gamma$ fixing infinity in $\HH^3$.  For any $g \in \Gamma \setminus
\Gamma_\infty$, the isometric sphere $S_g$ has boundary that is a
circle on the boundary $\CC$ at infinity of $\HH^3$.  This circle
bounds an open disk $D_g$ in $\CC$.  Two isometric spheres $S_g$ and
$S_h$ are \emph{visibly tangent} if their corresponding disks $D_g$
and $D_h$ are tangent on $\CC$, and for any other $k \in \Gamma
\setminus \Gamma_\infty$, the point of tangency is not contained in
the open disk $D_k$.

\begin{proof}
Suppose $\rho$ is not minimally parabolic.  Then it must have a rank 1
cusp.  Apply an isometry to $\HH^3$ so that the point at infinity
projects to this rank 1 cusp.  The Ford domain becomes a finite sided
region $P$ meeting this cusp.  Take a horosphere about infinity.
Because the Ford domain is finite sided, we may take this horosphere
about infinity sufficiently small that the intersection of the
horosphere with $P$ gives a subset of Euclidean space with sides
identified by elements of $\rho(\pi_1(C))$, conjugated appropriately.

The side identifications of this subset of Euclidean space, given by
the side identifications of $P$, generate the fundamental group of the
cusp.  But this is a rank 1 cusp, hence its fundamental group is
$\ZZ$.  Therefore, the side identification is given by a single
Euclidean translation.  The Ford domain $P$ intersects this horosphere
in an infinite strip, and the side identification glues the strip into
an annulus.  Note this implies two faces of $P$ are tangent at
infinity.

Now conjugate back to our usual view of $\HH^3$, with the point at
infinity projecting to the rank 2 cusp of the $(1,2)$--compression
body $\HH^3/\rho(\pi_1(C))$.  The two faces of $P$ tangent at infinity
are taken to two isometric spheres of the Ford domain, tangent at a
visible point on the boundary at infinity.
\end{proof}

\begin{remark}
The converse to Lemma \ref{lemma:min-parabolic} is not true.  There
exist examples of geometrically finite representations for which two
visible isometric spheres are visibly tangent, and yet the
representation is still minimally parabolic.  We see examples of this
in \cite{lackenby-purcell-2}.
\end{remark}

We next prove a result which will help us identify representations
which are \emph{not} discrete.

\begin{lemma}
Let $\Gamma$ be a discrete, torsion free subgroup of ${\rm PSL}(2,\CC)$ such
that $M = \HH^3/\Gamma$ has a rank two cusp.  Suppose that the point
at infinity projects to the cusp, and let $\Gamma_\infty$ be its
stabilizer in $\Gamma$.  Then for all $\zeta \in \Gamma \setminus
\Gamma_\infty$, the isometric sphere of $\zeta$ has radius at most the
minimal (Euclidean) translation length of all elements in
$\Gamma_\infty$.
\label{lemma:not-gf}
\end{lemma}

\begin{proof}
By the Margulis lemma, there exists an embedded horoball neighborhood
of the rank 2 cusp of $\HH^3/\Gamma$.  Let $H_\infty$ be a horoball
about infinity in $\HH^3$ that projects to this embedded horoball.
Let $\tau$ be the minimum (Euclidean) translation length of all
nontrivial elements in the group $\Gamma_\infty$, say $\tau$ is the
distance translated by the element $w_\tau$.  Suppose $S_\zeta$ has
radius $R$ strictly larger than $\tau$.  Without loss of generality,
we may assume $S_\zeta$ is visible, for otherwise there is some
visible face $S_\xi$ which covers the highest point of $S_\zeta$,
hence must have even larger radius.

Because the radius $R$ of $S_\zeta$ is larger than $\tau$, $S_\zeta$
must intersect $w_\tau (S_\zeta) = S_{\zeta w_\tau^{-1}}$, and in
fact, the center $w_\tau\zeta^{-1}(\infty)$ of $S_{\zeta w_\tau^{-1}}$
must lie within the boundary circle $S_\zeta \cap \CC$.

Consider the set of points $P$ equidistant from $\zeta^{-1}(H_\infty)$
and $w_\tau \zeta^{-1}(H_\infty)$.  Because these horoballs are the
same size, $P$ must be a vertical plane in $\HH^3$ which lies over the
perpendicular bisector of the line segment running from
$\zeta^{-1}(\infty)$ to $w_\tau\zeta^{-1}(\infty)$ on $\CC$.

Now apply $\zeta$.  This will take the plane $P$ to $S_0:=S_{\zeta
w_\tau^{-1} \zeta^{-1}}$.  We wish to determine the (Euclidean) radius
of $S_0$.  By Lemma \ref{lemma:isosphere-invert}, applying $\zeta$ is
the same as applying an inversion in $S_\zeta$, followed by a
Euclidean isometry.  Only the inversion will affect the radius of
$S_0$.  Additionally, the radius is independent of the location of the
center of the isometric sphere $S_\zeta$, so we may assume without
loss of generality that the center of $S_\zeta$ is at $0 \in \CC$ and
that the center of $S_{\zeta w_\tau^{-1}}$ is at $\tau \in \CC$.  Now
inversion in a circle of radius $R$ centered at zero takes the point
$\tau$ to $R^2/\tau$, and the point at infinity to $0$.  Thus the
center of $S_0$, which is the image of $\tau$ under $\zeta$, will be
of distance $R^2/\tau$ from a point on the boundary of $S_0$, i.e. the
image of $\infty$ on $P$ under $\zeta$.  Hence the radius of $S_0$ is
$R^2/\tau > R$.  Denote $R^2/\tau$ by $R_0$.  We have $R_0 > R >
\tau$.

Now we have a new face $S_0$ with radius $R_0>R>\tau$.  Again we may
assume it is visible.  The same argument as above implies there is
another sphere $S_1$ with radius $R_1>R_0>\tau$.  Continuing, we
obtain an infinite collection of visible faces of increasing radii.
These must all be distinct.  But this is impossible: an infinite
number of distinct faces of radius greater than $\tau$ cannot fit
inside a fundamental domain for $\Gamma_\infty$.  Thus $\Gamma$ is
indiscrete.
\end{proof}

The following lemma gives us a tool to identify the Ford domain of a
geometrically finite manifold.

\begin{lemma}
Let $\Gamma$ be a subgroup of ${\rm PSL}(2,\CC)$ with rank 2 subgroup
$\Gamma_\infty$ fixing the point at infinity.  Suppose the isometric
spheres corresponding to a finite set of elements of $\Gamma$, as well
as a vertical fundamental domain for $\Gamma_\infty$, cut out a
fundamental domain $P$ for $\Gamma$.  Then $\Gamma$ is discrete and
geometrically finite, and $P$ must be a Ford domain of $\HH^3/\Gamma$.
\label{lemma:finding-ford}
\end{lemma}

\begin{proof}
The discreteness of $\Gamma$ follows from Poincar{\'e}'s polyhedron
theorem.  The fact that it is geometrically finite follows directly
from the definition.

Suppose $P$ is not a Ford domain.  Since the Ford domain is only
well--defined up to choice of fundamental region for $\Gamma_\infty$,
there is a Ford domain $F$ with the same choice of vertical
fundamental domain for $\Gamma_\infty$ as for $P$.  Since $P$ is not a
Ford domain, $F$ and $P$ do not coincide.  Because both are cut out by
isometric spheres corresponding to elements of $\Gamma$, there must be
additional visible faces that cut out the domain $F$ than just those
that cut out the domain $P$.  Hence $F$ is a strict subset of $P$, and
there is some point $x$ in $\HH^3$ which lies in the interior of $P$,
but does not lie in the Ford domain.

But now consider the covering map $\phi\co \HH^3 \to \HH^3/\Gamma$.
This map $\phi$ glues both $P$ and $F$ into the manifold
$\HH^3/\Gamma$, since they are both fundamental domains for $\Gamma$.
So consider $\phi$ applied to $x$.  Because $x$ lies in the interior
of $P$, and $P$ is a fundamental domain, there is no other point of
$P$ mapped to $\phi(x)$.  On the other hand, $x$ does not lie in the
Ford domain $F$.  Thus there is some preimage $y$ of $\phi(x)$ under
$\phi$ which does lie in $F$.  But $F$ is a subset of $P$.  Hence we
have $y \neq x$ in $P$ such that $\phi(x) = \phi(y)$.  This
is a contradiction.
\end{proof}

\subsection{The {F}ord spine}

When we glue the Ford domain into the manifold $M=\HH^3/\Gamma$, as in
the proof of Lemma \ref{lemma:finding-ford}, the faces of the Ford
domain will be glued together in pairs to form $M$.  %% These glued faces
%% (and edges and vertices) form a geometric spine in $M$.

\begin{define}
The \emph{Ford spine} of $M$ is defined to be the image of the visible
faces of $\F$ under the covering $\HH^3\to M$.
\label{def:spine}
\end{define}

\begin{remark}
A spine usually refers to a subset of the manifold onto which there is
a retraction of the manifold.  Using that definition, the Ford spine
is not strictly a spine.  However, the Ford spine union the genus 2
boundary $\D_+C$ will be a spine for the compression body.
\end{remark}

Let $\rho$ be a geometrically finite uniformization.  Recall that the
\emph{domain of discontinuity} $\Omega_{\rho(\pi_1(C))}$ is the
complement of the limit set of $\rho(\pi_1(C))$ in the boundary at
infinity $\D_\infty \HH^3$.  See, for example, Marden \cite[section
2.4]{marden-book}.

\begin{lemma}
Let $\rho$ be a minimally parabolic geometrically finite
uniformization of a $(1,2)$--compression body $C$.  Then the manifold
$(\HH^3 \cup \Omega_{\rho(\pi_1(C))})/\rho(\pi_1(C))$ retracts onto
the boundary at infinity $(\bar{\F} \cap \CC)/ \Gamma_\infty$, union
the Ford spine.
\label{lemma:spine}
\end{lemma}

\begin{proof}
Let $H$ be a horosphere about infinity in $\HH^3$ that bounds a
horoball which projects to an embedded horoball neighborhood of the
cusp of $\HH^3/\rho(\pi_1(C))$.  Let $x$ be any point in $\F \cap
\HH^3$.  The nearest point on $H$ to $x$ lies on a vertical line
running from $x$ to infinity.  These vertical lines give a foliation
of $\F$.  All such lines have one endpoint on infinity, and the other
endpoint on $\bar{\F} \cap \CC$ or an isometric sphere of $\F$.
We obtain our retraction by mapping the point $x$ to the endpoint of
its associated vertical line, then quotienting out by the action of
$\rho(\pi_1(C))$.
\end{proof}

To any face $F_0$ of the Ford spine, we obtain an associated
collection of visible elements of $\Gamma$: those whose isometric
sphere projects to $F_0$ (or more carefully, a subset of their
isometric sphere projects to the face $F_0$).  We will often say that
an element $g$ of $\Gamma$ \emph{corresponds} to a face $F_0$ of the
Ford spine of $M$, meaning $S_g$ is visible, and (the visible subset
of) $S_g$ projects to $F_0$.  Note that if $g$ corresponds to $F_0$,
then so does $g^{-1}$ and $w_0 g^{\pm 1} w_1$ for any words $w_0, w_1
\in \Gamma_\infty$.

\section{Ford domains of compression bodies}\label{sec:ford}
%% \subsection{Ford domains of compression bodies}

Let $C$ be a $(1,2)$--compression body.  The fundamental group
$\pi_1(C)$ is isomorphic to $(\ZZ\times\ZZ)\ast \ZZ$.  The $\ZZ \times
\ZZ$ factor has generators $\alpha$ and $\beta$, and the generator of
the $\ZZ$ factor is $\gamma$.

Suppose $\rho\co \pi_1(C) \to {\rm PSL}(2,\CC)$ is a minimally parabolic
geometrically finite uniformization of $C$.  Then $\rho(\alpha)$ and
$\rho(\beta)$ are parabolic, and we will assume they fix the point at
infinity in $\HH^3$.  Together, they generate $\Gamma_\infty$.  The
third element, $\rho(\gamma)$, is a loxodromic element.  In
$\pi_1(C)$, $\alpha$ and $\beta$ are represented by loops in
$\partial_-C$.  To form the $(1,2)$--compression body, we add to
$\partial_-C \times I$ a 1--handle.  Then $\gamma$ is represented by a
loop around the core of this 1--handle.

%\begin{example}
In the simplest possible case imaginable, the Ford spine of
$\HH^3/\Gamma$ consists of a single face, corresponding to
$\rho(\gamma)$.  Note if this case happened to occur, then in the lift
to $\HH^3$, the only visible isometric spheres would correspond to
$\rho(\gamma)$, $\rho(\gamma^{-1})$, and their translates by elements
of $\Gamma_\infty$.  Cutting out regions bounded by these hemispheres
would give the region $\F$.  Topologically, the manifold
$\HH^3/\Gamma$ is obtained as follows.  First take $\F/\Gamma_\infty$.
The interior of $\F/\Gamma_\infty$ is homeomorphic to $T^2 \times
(0,\infty)$.  On the boundary on $\CC$ of $\F/\Gamma_\infty$ lie
two hemispheres, corresponding to $\rho(\gamma)$ and
$\rho(\gamma^{-1})$.  These are glued via $\rho(\gamma)$ to form
$\HH^3/\Gamma$ from $\F/\Gamma_\infty$.

\begin{figure}
\begin{center}
\begin{tabular}{ccc}
	\makebox{\begin{tabular}{c}
				\input{figures/simple-forddomain.pstex_t} \\
%%				\vspace{.1in} 
			\end{tabular}	} & 
	\hspace{1in} &
	\makebox{\begin{tabular}{c}
			\vspace{-0.2in} \\
			\includegraphics[width=1.5in]{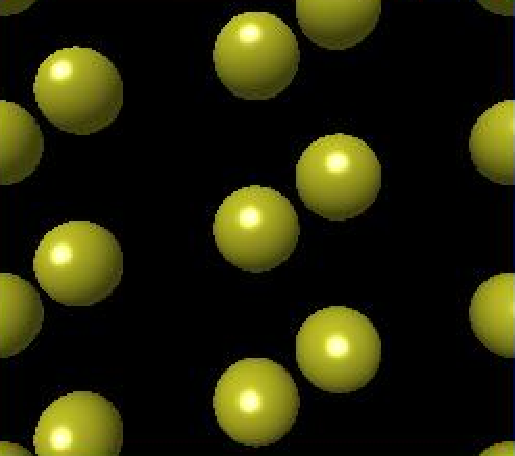}
			\end{tabular}}
\end{tabular}
\end{center}
\caption{Left: Schematic picture of a simple Ford domain.  Right:
	Three dimensional view of $\F$ in $\HH^3$.}
\label{fig:simple-ford}
\end{figure}

This situation is illustrated in Figure \ref{fig:simple-ford}. 
%\label{example:simple}
%\end{example}

%% NOTE:  STICKING IN ALL THESE BOXES AND TABULARS BELOW IS THE ONLY
%% WAY I CAN SEEM TO GET THE FORMATTING TO WORK THE WAY I WANT.  TRY
%% TO FIX LATER?

In the following lemma, we show that this simple Ford domain does, in
fact, occur.

\begin{lemma}
Let $C$ be a $(1,2)$--compression body.  There exists a minimally
parabolic geometrically finite uniformization of $C$, $\rho\co
\pi_1(C) \to {\rm PSL}(2,\CC)$ such that the Ford spine of
$\HH^3/\rho(\pi_1(C))$ consists of a single face, corresponding to the
loxodromic generator.
\label{lemma:simple-ford}
\end{lemma}

\begin{proof}
We construct such a structure by choosing $\rho(\alpha)$,
$\rho(\beta)$, $\rho(\gamma)$ in ${\rm PSL}(2,\CC)$.

Let $c \in \CC$ be such that $|c| > 2$, and let $\rho(\alpha)$,
$\rho(\beta)$, and $\rho(\gamma)$ be defined by
$$\rho(\alpha) = \mat{1&2|c|\\0&1}, \quad \rho(\beta) =
\mat{1&2i|c|\\0&1}, \quad \rho(\gamma)=\mat{c&-1\\1&0}.$$
Let $\Gamma$ be the subgroup of ${\rm PSL}(2, \CC)$ generated by
$\rho(\alpha)$, $\rho(\beta)$, and $\rho(\gamma)$.  By Lemma
\ref{lemma:iso-center-rad}, $S_{\rho(\gamma)}$ has center $0$, radius
$1$, and $S_{\rho(\gamma^{-1})}$ has center $c\in\CC$, radius $1$.
For $|c|>2$, $S_{\rho(\gamma)}$ will not meet $S_{\rho(\gamma^{-1})}$.
Note also that by choice of $\rho(\alpha)$, $\rho(\beta)$, all
translates of $S_{\rho(\gamma)}$ and $S_{\rho(\gamma^{-1})}$ under
$\Gamma_\infty$ are disjoint.  We claim that $\rho$ satisfies
the conclusions of the lemma.

Select a vertical fundamental domain for $\Gamma_\infty$ which
contains the isometric spheres $S_{\rho(\gamma)}$ and
$S_{\rho(\gamma^{-1})}$ in its interior.  This is possible by choice
of $\rho(\alpha)$, $\rho(\beta)$, and $\rho(\gamma)$.

Consider the region $P$ obtained by intersecting this fundamental
region with the complement of $B_{\rho(\gamma)}$ and
$B_{\rho(\gamma^{-1})}$.  As in the discussion above, we may glue this
region $P$ into a manifold $C_0$ by gluing $S_{\rho(\gamma)}$ to
$S_{\rho(\gamma^{-1})}$ via $\rho(\gamma)$, and by gluing vertical
faces by appropriate parabolic elements.  The manifold $C_0$ will be
homeomorphic to the interior of a $(1,2)$--compression body.

Then Poincar{\'e}'s polyhedron theorem implies that the manifold $C_0$
has fundamental group generated by $\rho(\alpha)$, $\rho(\beta)$, and
$\rho(\gamma)$.  Hence $C_0$ is the manifold $\HH^3/\Gamma$.

By Lemma \ref{lemma:finding-ford}, this fundamental region $P$ must
actually be the Ford domain for the manifold, and $\Gamma$ is
geometrically finite.  Since these isometric spheres are nowhere
tangent, $\rho$ is minimally parabolic, by Lemma
\ref{lemma:min-parabolic}.
\end{proof}

The examples of Ford domains that will interest us will be more
complicated than that in %Example \ref{example:simple}.
Lemma \ref{lemma:simple-ford}

\begin{example}
Fix $R>0$, and select $\varepsilon \in \RR$ so that
$0<\varepsilon<e^{-R}$, or equivalently, so that $\log(1/\varepsilon) \in
(R, \infty)$.  Set $\rho(\gamma)$ equal to
\begin{equation}
\mat{\displaystyle{\frac{i(1+\varepsilon)}{\sqrt{\varepsilon}}} &
	\displaystyle{ \frac{i}{\sqrt{\varepsilon}}}\vspace{.1in} \\
	-\displaystyle{\frac{i}{\sqrt{\varepsilon}}} &
	-\displaystyle{\frac{i}{\sqrt{\varepsilon}}}}.
\label{eqn:gamma-mat}
\end{equation}
Note that with $\rho(\gamma)$ defined in this manner, we have
$$\rho(\gamma^2) = \mat{-2-\varepsilon & -1 \\ 1 & 0}.$$
Thus the isometric sphere of $\rho(\gamma)$ has radius
$1/|i/\sqrt{\varepsilon}| = \sqrt{\varepsilon}$, while that of
$\rho(\gamma^2)$ has radius $1$ by Lemma \ref{lemma:iso-center-rad}.
Now select $\rho(\alpha)$ and $\rho(\beta)$ to be parabolic
translations fixing the point at infinity, with translation distance
large enough that the isometric spheres of $\rho(\gamma^2)$,
$\rho(\gamma^{-2})$, $\rho(\gamma)$, and $\rho(\gamma^{-1})$ do not
meet any of their translates under $\rho(\alpha)$ and $\rho(\beta)$.
The following will do:
$$\rho(\alpha) = \mat{1& 20 \\ 0 & 1}, \quad
\rho(\beta) = \mat{1 & 20i\\ 0& 1}.$$
\label{example:long-tunnel}
\end{example}

\begin{lemma}
The representation $\rho\co \pi_1(C) \to {\rm PSL}(2,\CC)$ defined in
Example \ref{example:long-tunnel} is a minimally parabolic
geometrically finite hyperbolic uniformization of %% the
%% $(1,2)$--compression body
$C$ whose Ford spine consists of exactly two
faces, corresponding to $\rho(\gamma)$ and $\rho(\gamma^2)$.
\label{lemma:gf}
\end{lemma}

\begin{proof}
Consider the isometric spheres corresponding to $\rho(\gamma)$,
$\rho(\gamma^{-1})$, $\rho(\gamma^2)$, and $\rho(\gamma^{-2})$.  We
will show that these faces, along with the faces of a vertical
fundamental domain for the action of $\rho(\alpha)$ and $\rho(\beta)$,
are the only faces of the Ford domain of the manifold
$\HH^3/\rho(\pi_1(C))$.  Since the faces corresponding to
$\rho(\gamma)$ and to $\rho(\gamma^{-1})$ glue together, and since the
faces corresponding to $\rho(\gamma^2)$ and to $\rho(\gamma^{-2})$
glue, the Ford domain glues to give a Ford spine with exactly two
faces.  The fact that the manifold is geometrically finite will then
follow by Lemma \ref{lemma:finding-ford}.

Choose vertical planes that cut out a vertical fundamental domain for
the action of $\Gamma_\infty$ and that avoid the isometric
spheres corresponding to $\rho(\gamma^{\pm 1})$ and $\rho(\gamma^{\pm
2})$.  Because the translation distances of $\rho(\alpha)$ and
$\rho(\beta)$ are large with respect to the radii of these isometric
spheres, this is possible.  For example, the planes $x=-10$, $x=10$,
$y=-10$, $y=10$ in $\{(x,y,z)|z>0\} = \HH^3$ will do.

Now, the isometric spheres of $\rho(\gamma)$ and $\rho(\gamma^{-1})$
have center $-1$ and $-1-\varepsilon$, respectively, and radius
$\sqrt{\varepsilon}$, by Lemma \ref{lemma:iso-center-rad}.  Similarly,
the isometric spheres of $\rho(\gamma^{2})$ and $\rho(\gamma^{-2})$
have centers $0$ and $-2-\varepsilon$, respectively, and radius $1$.
Then one may check: The isometric sphere of $\rho(\gamma^2)$ meets
that of $\rho(\gamma)$ in the plane $x=-1+\varepsilon/2$.  The isometric
sphere of $\rho(\gamma)$ meets that of $\rho(\gamma^{-1})$ in the
plane $x=-1-\varepsilon/2$, and the isometric sphere of
$\rho(\gamma^{-1})$ meets that of $\rho(\gamma^{-2})$ in the plane
$x=-1-3\varepsilon/2$, as in Figure \ref{fig:dumbbell}.  These are the
only intersections of these spheres that are visible from infinity.
If we glue the isometric spheres of $\rho(\gamma^{\pm 1})$ via
$\rho(\gamma)$ and the isometric spheres of $\rho(\gamma^{\pm 2})$ via
$\rho(\gamma^2)$, then these three edges of intersection are all glued
to a single edge.  %% One may
%% check that the angles around this edge sum to $2\pi$.

\begin{figure}
\begin{center}
	\input{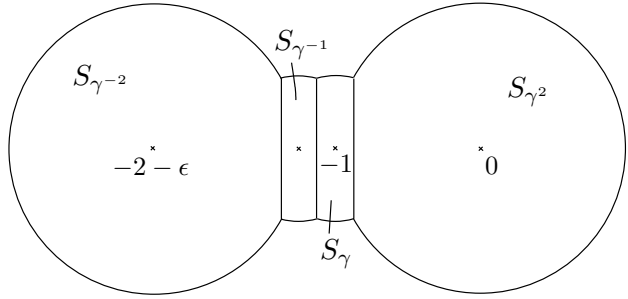}
\end{center}
\caption{The isometric spheres corresponding to $\rho(\gamma^{-2})$,
	$\rho(\gamma^{-1})$, $\rho(\gamma)$, and $\rho(\gamma^2)$.}
\label{fig:dumbbell}
\end{figure}

Consider the monodromy around this edge.  We must show that it is the
identity.  Note that a meridian of the edge is divided into three
arcs, running from the faces labeled $S_{\gamma^{-1}}$ to
$S_{\gamma^{-2}}$, from $S_{\gamma^2}$ to $S_{\gamma}$, and from
$S_{\gamma^{-1}}$ to $S_{\gamma}$.  To patch the first pair of arcs
together, we glue $S_{\gamma^{-2}}$ to $S_{\gamma^2}$ using the
isometry $\gamma^{-2}$.  To patch the second and third pairs of arcs,
we glue $S_{\gamma}$ to $S_{\gamma^{-1}}$ by the isometry $\gamma$.
The composition of these three isometries is
$\gamma^{-2}\gamma\gamma$, which is the identity, as required.

Hence, by Poincar\'e's polyhedron theorem, the space obtained by
gluing faces of the polyhedron $P$ cut out by the above isometric
spheres and vertical planes is a manifold, with fundamental group
generated by $\rho(\gamma)$, $\rho(\gamma^2)$, $\rho(\alpha)$, and
$\rho(\beta)$.

We need to show that this is a uniformization of $C$, i.e., that
$\HH^3/\rho(\pi_1(C))$ is homeomorphic to the interior of $C$.  The
Ford spine of $\HH^3/\rho(\pi_1(C))$ has two faces, one of which has
boundary which is the union of the 1--cell of the spine and an arc on
$\partial_+C$.  Collapse the 1--cell and this face.  The result is a
new complex with the same regular neighborhood.  It now has a single
2--cell attached to $\partial_+C$.  Thus, $\HH^3/\rho(\pi_1(C))$ is
obtained by attaching a 2--handle to $\partial_+C \times I$, and then
removing the boundary.  In other words, $\HH^3/\rho(\pi_1(C))$ is
homeomorphic to the interior of $C$. 

Thus $\HH^3/\rho(\pi_1(C))$ is homeomorphic to the interior of $C$,
and has a convex fundamental domain $P$ with finitely many faces.  By
Lemma \ref{lemma:finding-ford}, this convex fundamental domain $P$ is
actually the Ford domain.  Finally, since none of the isometric
spheres of the Ford domain are visibly tangent at their boundaries, by
Lemma \ref{lemma:min-parabolic} the representation is minimally
parabolic.  Hence it is a minimally parabolic geometrically finite
uniformization of $C$.
\end{proof}

\subsection{Dual edges}

To each face of the Ford spine, there is an associated dual edge,
which is defined as follows.  For any face $F$ of the Ford spine,
there is some $g \in \Gamma \backslash \Gamma_\infty$ such that
(subsets of) $S_g$ and $S_{g^{-1}}$ are faces of a Ford domain, and
$S_g$ and $S_{g^{-1}}$ project to $F$.  Above each of these isometric
spheres lies a vertical arc, running from the top of the isometric
sphere (i.e. the geometric center of the hemisphere) to the point at
infinity.  Define the dual edge to be the union of the image of these
two arcs in $\HH^3/\rho(\pi_1(C))$.

\begin{figure}
\begin{center}
	\input{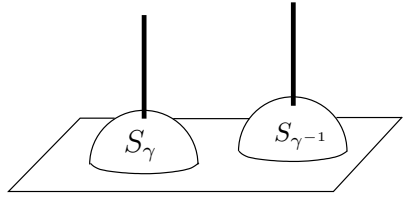}
\end{center}
\caption{The dual to the simplest Ford spine is an edge that lifts to
	a collection of vertical geodesics in $\F$, shown in bold.}
\label{fig:simple-dual}
\end{figure}

\begin{lemma}
For any uniformization $\rho\co \pi_1(C) \to {\rm PSL}(2,\CC)$, the
core tunnel will be homotopic to the edge dual to the isometric sphere
corresponding to the loxodromic generator of $\rho(\pi_1(C))$.
\label{lemma:homotopic-unknotting}
\end{lemma}

\begin{proof}
Denote the loxodromic generator by $\rho(\gamma)$.  Consider the core
tunnel in the compression body $\HH^3/\rho(\pi_1(C))$.  Take a
horoball neighborhood $H$ of the cusp, and its horospherical torus
boundary.  The core tunnel runs through this horospherical torus $\D
H$, into the cusp.  Denote by $\tilde{H}$ a lift of $H$ to $\HH^3$
about the point at infinity in $\HH^3$.

There is a homeomorphism from $C$ to $\HH^3/\rho(\pi_1(C)) \setminus
H$.  Slide the tunnel in $C$ so that it starts and ends at the same
point, and so that the resulting loop represents $\gamma$.  The image
of this loop under the homeomorphism to $\HH^3/\rho(\pi_1(C))
\setminus H$ is some loop.  This lifts to an arc in $\HH^3$ starting
on $\tilde{H}$ and ending on $\rho(\gamma)(\tilde{H})$.  Extend this
to an arc in $\HH^3 / \rho(\pi_1(C))$ by attaching a geodesic in
$\tilde{H}$ and in $\rho(\gamma)(\tilde{H})$.  This is isotopic to
(the interior of) the core tunnel.  Now homotope this to a geodesic.
It will run through the isometric sphere corresponding to
$\rho(\gamma^{-1})$ once.
\end{proof}

\section{Long unknotting tunnels}\label{sec:long-tunnels}
We are now ready to give the geometric proof of our main theorem.

\begin{theorem}
There exist finite volume one--cusped hyperbolic tunnel number one
manifolds for which the geodesic representative of the unknotting
tunnel is arbitrarily long, as measured between the maximal horoball
neighborhood of the cusp.
\label{thm:long}
\end{theorem}

Recall that a \emph{tunnel number one manifold} is a manifold $M$ with
torus boundary components which admits an unknotting tunnel, that is,
a properly embedded arc $\tau$, the exterior of which is a handlebody.

Recall also that the length of the geodesic representative of an
unknotting tunnel is measured outside a maximal horoball neighborhood
of the cusp.

Before proving Theorem \ref{thm:long}, we need to prove a similar
statement for minimally parabolic geometrically finite hyperbolic
uniformizations of a $(1,2)$--compression body.

%% FOR THIS PROOF, WE NEED:
%% Fundamental group of compression body and representations
%%  i.e. define alpha, beta, gamma
%% Poincare's theorem giving geometrically finite structure
%% Isometric spheres and length of unknotting tunnel
%% Relationship of unknotting tunnel and isometric sphere of $\gamma$
%% Visible intersections and manifolds

\begin{prop}
	For any $R>0$, there exists a minimally parabolic geometrically
	finite uniformization of a $(1,2)$--compression body such that the
	geodesic representative of the homotopy class of the core tunnel has
	length at least $R$.
\label{prop:geom-fin-long}
\end{prop}

\begin{proof}
We will prove Proposition \ref{prop:geom-fin-long} by finding an
explicit minimally parabolic geometrically finite uniformization of a
$(1,2)$--compression body $C$.  For fixed $R>0$, our explicit
uniformization will be that given in Example \ref{example:long-tunnel}
above.  By Lemma \ref{lemma:gf}, this is a minimally parabolic
geometrically finite hyperbolic uniformization of the
$(1,2)$--compression body $C$ whose Ford spine consists of exactly two
faces, corresponding to $\rho(\gamma)$ and $\rho(\gamma^2)$.  We claim
that the geodesic representative of the homotopy class of the core
tunnel has length at least $R$.

\begin{lemma}
Let $\rho\co \pi_1(C) \to {\rm PSL}(2,\CC)$ be a discrete, faithful
representation such that $\rho(\alpha)$, $\rho(\beta)$ are parabolics
fixing the point at infinity in $\HH^3$, and $\rho(\gamma)$ is as in
equation (\ref{eqn:gamma-mat}).  Then the geodesic representative of
the homotopy class of the core tunnel has length greater than
$R$.
\label{lemma:gf-long}
\end{lemma}

\begin{proof}
By Lemma \ref{lemma:homotopic-unknotting}, the core tunnel is
homotopic to the geodesic dual to the isometric spheres corresponding
to $\rho(\gamma)$ and $\rho(\gamma^{-1})$.  The length of this
geodesic is twice the distance along the vertical geodesic from the
top of one of the isometric spheres corresponding to $\rho(\gamma^{\pm
1})$ to a maximal horoball neighborhood of the cusp about infinity.
Since the isometric sphere of $\rho(\gamma^2)$ has radius $1$, a
maximal horoball about the cusp will have height at least $1$.  The
isometric sphere of $\rho(\gamma)$ has radius $\sqrt{\varepsilon}$.
Integrating $1/z$ from $z=\sqrt{\varepsilon}$ to $1$, we find that the
distance along this vertical arc is at least
$\log{1/\sqrt{\varepsilon}}$.  Hence the length of the geodesic
representative of the core tunnel is at least $\log{1/\varepsilon}$.
By choice of $\varepsilon$, this length is greater than $R$.
\end{proof}

\begin{remark}
Note in the proof above that we may strengthen Lemma
\ref{lemma:gf-long} as follows.  Because of the choice of $\varepsilon$
in equation (\ref{eqn:gamma-mat}), there exists some neighborhood $U$
of the matrix of (\ref{eqn:gamma-mat}) such that if $\rho(\alpha)$,
$\rho(\beta)$ are as above, but $\rho(\gamma)$ lies in $U$, then the
geodesic representative of the homotopy class of the core tunnel has
length greater than $R$.
\end{remark}

This completes the proof of Proposition \ref{prop:geom-fin-long}.
\end{proof}

%% We don't need to show there exists a geometrically infinite
%% structure.  All we really need is that maximal cusps are dense in
%% the boundary, and that we can get to the boundary.

Before we present the proof of Theorem \ref{thm:long}, we need to
recall terminology from Kleinian group theory.

We define the (restricted) character variety $V(C)$ to be the space of
conjugacy classes of representations $\rho \co \pi_1(C) \to {\rm
PSL}(2, \CC)$ such that elements of $\pi_1(\partial_- C)$ are sent to
parabolics.  Note this definition agrees with Marden's definition of
the representation variety in \cite{marden-book}, but is a restriction
of the character variety in Culler and Shalen's classic paper
\cite{culler-shalen}.  Convergence in $V(C)$ is known as algebraic
convergence.

Let $GF_0(C)$ denote the subset of $V(C)$ consisting of conjugacy
classes of minimally parabolic geometrically finite uniformizations of
$C$, given the algebraic topology.  It follows from work of Marden
\cite[Theorem 10.1]{marden:geom} that $GF_0(C)$ is an open subset of
$V(C)$.  We are interested in a type of structure that lies on the
boundary of $GF_0(C)$.  These structures are discrete, faithful
representations of $C$ that are geometrically finite, but not
minimally parabolic.

\begin{define}
A \emph{maximal cusp for $C$} is a geometrically finite uniformization
of $C$, $\rho\co \pi_1(C) \to {\rm PSL}(2,\CC)$ such that every component of
the boundary of the convex core of $\HH^3/\rho(\pi_1(C))$ is a
3--punctured sphere.
\end{define}

A maximal cusp is in some sense the opposite of a minimally parabolic
representation.  In a minimally parabolic representation, no elements
of $\D_+C$ are pinched.  In a maximal cusp, a full pants decomposition
of $\D_+C$, or the maximal number of elements, is pinched to parabolic
elements.

Due to a theorem of Canary, Culler, Hersonsky, and Shalen
\cite[Corollary 16.4]{cchs}, conjugacy classes of maximal cusps for
$C$ are dense on the boundary of $GF_0(C)$ in $V(C)$.  This theorem,
an extension of work of McMullen \cite{mcmullen:max-cusps}, is key in
the proof of Theorem \ref{thm:long}.

\begin{proof}[Proof of Theorem \ref{thm:long}]
Let $\rho_0$ be the geometrically finite representation of the proof
of Proposition \ref{prop:geom-fin-long}, with core tunnel homotopic to
a geodesic of length strictly greater than $R$.  The translation
lengths of $\rho_0(\alpha)$ and $\rho_0(\beta)$ are bounded, say by
$B$.

We will consider $\rho_0$ to be an element of the character variety
$V(C)$.  Indeed, define $\R$ to be the set of all representations
where $\rho(\alpha)$ and $\rho(\beta)$ are parabolics fixing infinity
with length bounded by $B$, and with $\rho(\gamma)$ fixed as in
equation (\ref{eqn:gamma-mat}).  If we view the character variety
$V(C)$ as a subset of the variety of representations $\rho$ of
$\pi_1(C)$ where $\rho(\alpha)$ and $\rho(\beta)$ have been suitably
normalized to avoid conjugation, then we may consider $\R$ as a subset
of $V(C)$.  Note $\rho_0$ is in $\R$.

The set $\R$ is clearly path connected.  By Lemma \ref{lemma:gf-long},
for all uniformizations of $C$ in $\R$, the length of the geodesic
representative of the core tunnel is at least $R$.

Moreover, notice that $\R$ includes indiscrete representations, as
follows.  Recall that the isometric sphere corresponding to $\gamma^2$
has radius $1$ when $\rho(\gamma)$ is defined as in equation
(\ref{eqn:gamma-mat}).  Thus by Lemma \ref{lemma:not-gf}, whenever the
translation length of $\alpha$ is less than $1$, the representation
cannot be discrete.

Then consider a path in $\R$ from $\rho_0$ to an indiscrete
representation.  At some point along this path, we come to $\R \cap \D
GF_0(C)$.  

By work of Canary, Culler, Hersonsky, and Shalen \cite{cchs},
generalizing work of McMullen \cite{mcmullen:max-cusps}, the set of
maximal cusps is dense in the boundary of geometrically finite
structures $\D GF_0(C)$.  

It follows that we can find a sequence of geometrically finite
representations $\rho_n$ of $\pi_1(C)$ such that the conformal
boundaries of the manifolds $C_n := \HH^3/\rho_n(\pi_1(C))$ are
maximally cusped genus two surfaces, $C_n$ are homeomorphic to the
interior of $C$, and such that the algebraic limit of these manifolds
$C_n$ is a manifold $M = \HH^3/\rho_\infty(\pi_1(C))$ where
$\rho_\infty$ is in $\R$.  By the remark following Lemma
\ref{lemma:gf-long}, for large $n$, the core tunnels of the $C_n$ will
have geodesic representative with length greater than $R$.

Now, there exists a maximally cusped hyperbolic structure on the genus
2 handlebody $H$.  In fact, by work of Canary, Culler, Hersonsky, and
Shalen \cite[Corollary 15.1]{cchs}, such structures are dense in the
boundary of geometrically finite structures on handlebodies.  Thus,
there exists a hyperbolic manifold $\HH^3/\Gamma_1$ homeomorphic to
the interior of $H$, such that every component of the boundary of the
convex core of $\HH^3/\Gamma_1$ is a 3--punctured sphere.  We will
continue to denote the hyperbolic manifold $\HH^3/\Gamma_1$ by $H$.

\begin{figure}
\begin{center}
\includegraphics[height=1.5in]{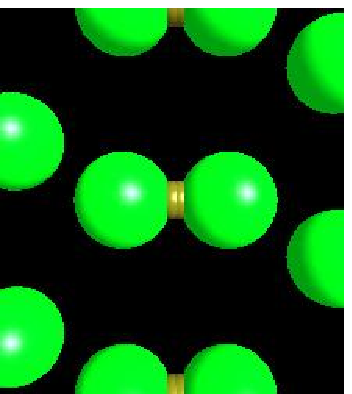} \hspace{.1in}
\includegraphics[height=1.5in]{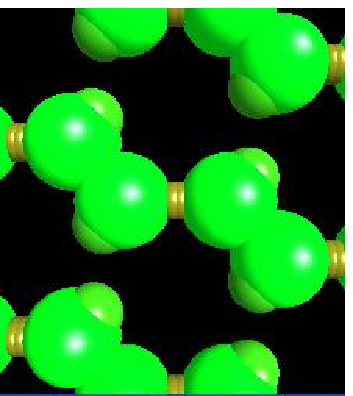} \hspace{.1in}
\includegraphics[height=1.5in]{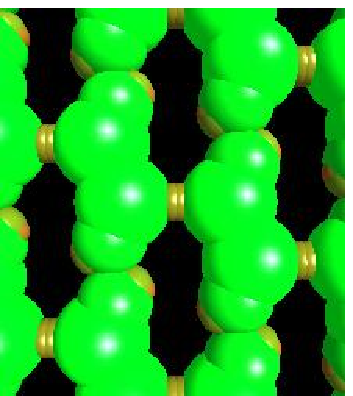} \hspace{.1in}
\includegraphics[height=1.5in]{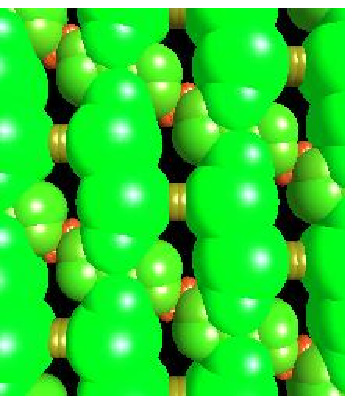}
\end{center}
\caption{Shown is a picture of $\F$ for four geometrically finite
	structures with long unknotting tunnel.  These structures are
	converging to a structure on $\D GF_0(\pi_1(C))$.  Note that in
	each of the four structures, the pattern of isometric spheres
	corresponding to that of Figure \ref{fig:dumbbell} is visible,
	although the number of visible isometric spheres increases.}
\label{fig:long-inf}
\end{figure}

Let $\phi_n$ be any homeomorphism from $\D_+C$ to $H$ taking
parabolics of $C_n$ on $\D_+C$ to the parabolics of $\D H$.  Because
$\phi_n$ takes 3--punctured spheres to 3--punctured spheres, it
extends to an isometry.  Hence we may glue $C_n$ to $H$ via $\phi_n$
and obtain a tunnel number one manifold with three drilled out curves,
corresponding to the three parabolics of $\D_+C$.  These are three
torus boundary components of $M_n:=C_n \cup_{\phi_n} H$. 

Select Dehn filling slopes $s^1$, $s^2$, $s^3$ on these three boundary
components that act as gluing one boundary to the other by a high
power of a Dehn twist.  When we Dehn fill along these slopes, the
result is a tunnel number one manifold $M_n(s^1, s^2, s^3)$.  By work
of Thurston \cite{thurston}, as the length of the slopes increases,
the Dehn filled manifold approaches $M_n$ in the geometric topology.
Thus the length of the geodesic representative of the homotopy class
of the unknotting tunnel in $M_n(s^1, s^2, s^3)$ approaches the length
of the geodesic representative of the homotopy class of the core
tunnel in $C_n$ as the lengths of $s^1$, $s^2$, and $s^3$ increase in
$M_n$.

Hence for large enough $n$ and long enough slopes $s^1$, $s^2$, $s^3$,
the Dehn filled manifold $M_n(s^1, s^2, s^3)$ is a tunnel number one
manifold with unknotting tunnel homotopic to a geodesic of length at
least $R$.
\end{proof}

\subsection{Remarks}
While Theorem \ref{thm:long} gives us a manifold whose unknotting
tunnel has a long geodesic representative, the proof does not
guarantee that this tunnel is isotopic to a geodesic, even if we could
guarantee that the core tunnel is isotopic to a geodesic in the
approximating geometrically finite structures $C_n$.  This isn't
important for the proof of Theorem \ref{thm:long}.  However, in
\cite{lackenby-purcell-2}, we will explain how to modify the above
proof so that the unknotting tunnel is isotopic to a geodesic.

\section{Knots in homology 3-spheres}\label{sec:homology}

In this section, we refine the construction in Theorem \ref{thm:long}
in order to control the homology of the resulting manifolds.

\begin{theorem}
There exist hyperbolic knots in homology 3-spheres which have
tunnel number one, for which the geodesic representative
of the unknotting tunnel is arbitrarily long.
\label{thm:homology}
\end{theorem}

The manifolds in the proof of Theorem \ref{thm:long} were contructed
by starting with maximally cusped geometrically finite uniformizations
of the compression body $C$ and the handlebody $H$, gluing them via an
isometry, and then performing certain Dehn fillings. We will now vary
this construction a little. We will again use maximally cusped
geometrically finite uniformizations of the $(1,2)$-compression body
$C$ and the genus 2 handlebody $H$, but we will not glue them
directly.  Instead, we will also find a maximally cusped geometrically
finite uniformization of $S \times I$, where $S$ is the closed
orientable surface with genus 2, and we will glue $C$ to $S \times \{
1 \}$ and glue $H$ to $S \times \{ 0 \}$.  In both gluings, the
parabolic loci will be required to match together, although we will
leave these loci unglued. The result is therefore a tunnel number one
manifold, with 6 disjoint embedded simple closed curves removed.  We
will then perform certain Dehn fillings along these 6 curves to give
the required tunnel number one manifolds. The choice of hyperbolic
structure on $H$ requires some care. In particular we will need the
following terminology and results.

Let ${\rm ML}(\partial H)$ (respectively, ${\rm PML}(\partial H)$) be
the space of measured laminations (respectively, projective measured
laminations) on $\partial H$. (See for example
\cite{fathi-laudenbach-poenaru}). Let $i( \cdot, \cdot)$ denote the
intersection number between measured laminations. A measured
lamination $\lambda$ is said to be {\sl doubly incompressible} if
there is an $\epsilon > 0$ such that $i(\lambda, \partial E) >
\epsilon$ for all essential annuli and discs $E$ properly embedded in
$H$. Similarly, a projective measured lamination is {\sl doubly
incompressible} if any of its inverse images in ${\rm ML}(\partial H)$
is doubly incompressible. It is a consequence of Thurston's
geometrization theorem \cite{morgan-smith} that if $P$ is a collection
of simple closed curves on $\partial H$ that are pairwise non-parallel
in $\partial H$, essential in $\partial H$ and doubly incompressible,
then there is a geometrically finite uniformization of $H$.  Let $P$
be the part of its parabolic locus $P$ that lies on $\partial_+C$. The
set of doubly incompressible projective measured laminations forms a
non-empty open subset of ${\rm PML}(\partial H)$ (see \cite{Lecuire}).

\begin{lemma}
\label{homtrivpa}
There is a homeomorphism
$\psi \colon \partial H \rightarrow \partial H$ satisfying
the following conditions:
\begin{enumerate}
\item $\psi$ is pseudo-Anosov;
\item its stable and unstable laminations are doubly incompressible;
\item the induced homomorphism $\psi_\ast \colon H_1(\partial H) \rightarrow
H_1(\partial H)$ is the identity.
\end{enumerate}
\end{lemma}

\begin{proof} 
Since the stable laminations of pseudo-Anosovs are dense in ${\rm
PML}(\partial H)$, and the set of doubly incompressible laminations is
open and non-empty, there is a pseudo-Anosov homeomorphism $g$ with
doubly incompressible stable lamination.  Let $h$ be a pseudo-Anosov
on $\partial H$ that acts trivially on $H_1(\partial H)$ (see
\cite{thurston-surfaces}).  Let $\lambda_+$ and $\lambda_-$ be its
stable and unstable projective measured laminations, which we may
assume are distinct from the unstable lamination of $g$. Then the
pseudo-Anosov $g^m h g^{-m}$ also acts trivially on $H_1(\partial H)$.
Its stable and unstable laminations are $g^m(\lambda_+)$ and
$g^m(\lambda_-)$, which are arbitrarily close to the stable lamination
of $g$ for large $m$. Hence, they too are doubly incompressible when
$m$ is large.  Thus, we may set $\psi$ to be one such $g^m h g^{-m}$.
\end{proof}

\begin{proof}[Proof of Theorem \ref{thm:homology}]
Let $\phi \colon \partial_+ C \rightarrow \partial H$ be a
homeomorphism such that, when $C$ is glued to $H$ via $\phi$, the
result is the standard genus two Heegaard splitting of the solid
torus.  Fix a maximally cusped geometrically finite uniformization of
$C$ from the proof of Theorem \ref{thm:long}, for which the core
tunnel has long geodesic representative. Let $P$ be its parabolic
locus.  Then $\phi(P)$ is a collection of simple closed curves on $H$.

Let $\tau$ be a composition of Dehn twists, one twist around each
component of $\phi(P)$. Let $\psi$ be the pseudo-Anosov homeomorphism
provided by Lemma \ref{homtrivpa}. By replacing $\psi$ by a power of
$\psi$ if necessary, we may assume that for each core curve $\alpha$
of $P$, $i(\alpha, \psi(\alpha)) \not= 0$.  The tunnel number one
manifold that we are aiming for is obtained by gluing $C$ to $H$ via
$\psi^{m} \tau^{-M} \psi^{-1} \tau^M \phi$ for large integers $m$ and
$M$.  Since $\psi$ acts trivially on homology, this has the same
homology as if we had glued by $\phi$, which gives the solid torus.
Thus, this manifold is indeed the exterior of a knot in a homology
3-sphere.

We first choose the integer $m$. As $m$ tends to infinity, $\psi^m
\phi(P)$ tends to the stable lamination of $\psi$ in ${\rm
PML}(\partial H)$. Hence, we may choose such an $m$ so that $\psi^m
\phi(P)$ is doubly incompressible.

We start with three manifolds:
\begin{enumerate}
\item $C - P$;
\item $(S \times [0,1]) - ((\phi (P) \times \{ 1 \}) \cup (\psi\phi (P) \times \{ 0 \}))$;
\item $H - \psi^{m} \phi(P)$.
\end{enumerate}
Here, $S$ is the genus two surface, which we identify with $\partial H$.
The second of the above manifolds has a geometrically finite
uniformization, by Thurston's geometrization theorem. This 
because any essential annulus in $S \times [0,1]$ with boundary in
$(\phi( P) \times \{ 1 \}) \cup (\psi\phi (P)  \times \{ 0 \})$ can be
homotoped, relative to its boundary, so that it lies entirely
in $(\phi (P) \times \{ 1 \}) \cup (\psi \phi(P) \times \{ 0 \})$.
Similarly, because $\psi^m \phi(P)$ is doubly incompressible,
$H - \psi^m \phi(P)$ admits a geometrically finite hyperbolic
structure. Glue $C - P$ to $(S - \phi(P)) \times \{1 \}$ via
$\phi$, and glue $(S - \psi\phi(P)) \times \{ 0 \}$ to
$H - \psi^m \phi(P)$ via $\psi^{m-1}$. Since these manifolds
have conformal boundary that consists of 3-punctured spheres,
this gluing can be performed isometrically. 

As in the proof of Theorem \ref{thm:long}, we now perform certain Dehn
fillings on the toral cusps of this manifold, apart from the cusp
corresponding to $\partial_- C$.  If the Dehn filling is done
correctly, this has the effect of modifying the gluing map by powers
of Dehn twists.  We may apply any iterate of these Dehn twists, and so
we apply the $M$th iterate, where $M$ is some large positive integer,
along each of the curves $\phi(P) \times \{ 1 \}$ in $S \times \{1\}$
and the $(-M)$th power along each of the curves in $\psi\phi(P) \times
\{ 0 \}$ in $S \times \{0 \}$. Thus, the gluing map becomes $\psi^{m}
\tau^{-M} \psi^{-1} \tau^M \phi$.  As $M$ tends to infinity, these
manifolds tend geometrically to the unfilled manifold. In particular,
for large $M$, the geodesic representative of its unknotting tunnel
will be long.
\end{proof}

\section{The {D}ehn filling construction}\label{sec:dehn}
In this section, we give the proof of Theorem \ref{thm:long} that uses
Dehn filling and homology.  

Let $X$ be a compact $3$--manifold with four torus boundary components
and of Heegaard genus $2$. This means there is a closed genus $2$
surface $F$ in the interior of $X$ which separates $X$ into two
compression bodies, each homeomorphic to the manifold $V$ obtained by
adding one $2$--handle onto a copy of $F\times[0,1]$ along an
essential separating simple closed curve in $F\times \{1\}.$ We label
the torus boundary components of $X$ by $A_0$, $A_1$, $B_0$, $B_1$ so
that $A_0$ and $B_0$ are on the same side of $F$.

Let $\beta_0$, $\beta_1$, and $\alpha_1$ be essential simple closed
curves on $B_0$, $B_1$, and $A_1$, respectively.  Let $M =
X(\alpha_1,\beta_0,\beta_1)$ be the manifold obtained by Dehn filling
using these slopes, so that $M$ has a single boundary component $A_0.$
Gluing a solid torus to each of the two boundary components of $V$
yields a genus $2$ handlebody. It follows that $M$ has tunnel number
one; indeed a tunnel is obtained using an arc with endpoints on $A_0$
that goes round the core of the solid torus used to fill along $B_0$.

\begin{lemma}
There exists $X$ as above such that the interior of $X$ admits a
complete hyperbolic structure of finite volume, such that $H_1(X)
\cong \Gamma_A \oplus \Gamma_B$ where $\Gamma_A \cong \Gamma_B \cong
\ZZ^2$, and under maps induced by inclusion, $H_1(A_i) = \Gamma_A$ and
$H_1(B_i) = \Gamma_B$ for $i=1,2.$
\label{lemma:hyp-link}
\end{lemma}

\begin{figure}
\input{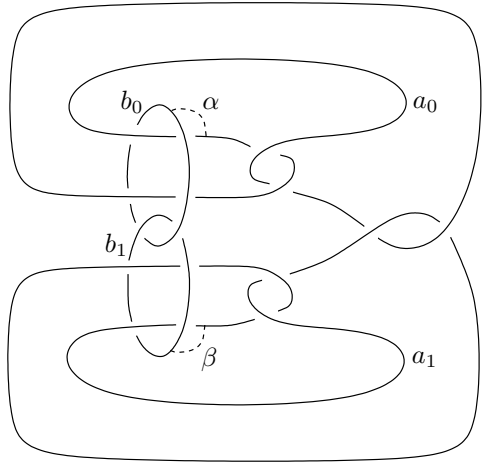}
\caption{A hyperbolic link satisfying the conditions of Lemma
	\ref{lemma:hyp-link}.} 
\label{fig:tunnellink}
\end{figure}

\begin{proof}
An example of $X$ is provided by the exterior of the $4$ component
link $L$ in $S^3$ shown in Figure \ref{fig:tunnellink}.  The link $L =
a_0 \cup a_1 \cup b_0 \cup b_1$ consists of two linked copies of the
Whitehead link and is hyperbolic (by SnapPea \cite{snappea}).
Furthermore $Lk(a_0,a_1) = Lk(b_0,b_1) = 1$ and $Lk(a_i,b_j)=0$.  The
diagram also shows disjoint arcs $\alpha$ connecting $a_0$ to $b_0$
and $\beta$ connecting $a_1$ to $b_1$.  It is easy to slide these arcs
and links in such a way that the pair of graphs $a_0 \cup b_0 \cup
\alpha$ and $a_1 \cup b_1 \cup \beta$ are spines of the handlebodies
of the genus $2$ Heegaard spliting of $S^3$.  It follows that $X =
S^3-{\eta}(L)$ has the required properties and $A_i =
\partial\eta(a_i)$ and $B_i = \partial\eta(b_i)$.  Here $\eta(L)$
denotes an open tubular neighborhood of $L$.
\end{proof}

Suppose $X$ is a link of Lemma \ref{lemma:hyp-link}.  Let $x$ be a
basepoint for $X$ on the boundary of a maximal horoball neighborhood
of the cusp corresponding to $A_0$.  The idea for finding the Dehn
fillings to give $M$ is that given a base point $x$ in $X$ and $R>0$,
there are only finitely many homotopy classes of loops in $X$ based at
$x$ with length at most $3R$.  These give finitely many classes in
$H_1(X)$ and hence, under projection, finitely many classes
$\gamma_1,\cdots \gamma_p \in \Gamma_B$.  The Dehn fillings used to
obtain $M$ are chosen so that $H_1(M) \cong \ZZ \oplus \ZZ_n$ with the
image of $\Gamma_B$ being $\ZZ_n$ and that of $\Gamma_A$ being $\ZZ$.
The fillings are also chosen so that none of the images of $\gamma_i$
generate $\ZZ_n$, and so that the hyperbolic metric on $M$ is
geometrically close to that of $X$.  In particular, we may assume that
there is a bilipschitz homeomorphism, with bilipschitz constant very
close to $1$, between the $R$ neighborhood of the basepoint of $x$ in
$X$ and a subset of $M$.  Let $m$ be the image of $x$ in $M$, which we
may assume lies on the boundary of a maximal horoball neighborhood of
the cusp of $M$.  Then every loop in $M$ based at $m$ of length at
most $2R$ corresponds to a loop in $X$ based at $x$ with length at
most $3R$, say.

\begin{lemma}
Suppose $\Gamma$ is a free abelian group of rank $2$ and $\gamma_1,
\cdots \gamma_p \in \Gamma$. Then there is an integer $n>0$ and an
epimorphism $\phi_n\co \Gamma \rightarrow \ZZ_n$ such that for all
$i$ the element $\phi_n(\gamma_i)$ does not generate $\ZZ_n$.
\label{lemma:not-generate}
\end{lemma}

\begin{proof}
Clearly we may assume that for all $i$, $\gamma_i \ne 0$.  Then we may
identify $\Gamma$ with $\ZZ^2$ so that $\gamma_i = (a_i, b_i)$ and for
all $i$, $a_i \ne 0$.  Set $m = \max_i |b_i|+2$ and define a
homomorphism $\phi\co \ZZ^2 \rightarrow \ZZ$ by $\phi((a,b)) = 2ma -
b$ which is surjective because $\phi((1,2m-1)) = 1$.  Set $c_i =
|\phi(\gamma_i)|$. Then
$$c_i = |2ma_i-b_i| \ge 2m|a_i|-|b_i| \ge 2m-m\ \ge 2,$$
using that $|a_i|\ge 1$ and $|b_i|\le m$ and $m\ge 2$.  Now define $n
= \prod_i c_i$ and define $\phi_n(\gamma) = \phi(\gamma)$ mod $n$.
Then $\phi_n(\gamma_i) = \pm c_i$ and $c_i \ne 1$ divides $n$ and
therefore does not generate $\ZZ_n.$
\end{proof}

\medskip
For the Dehn fillings of $B_0$ and $B_1$, choose simple closed curves
$\beta_i\subset B_i$ which generate the kernel of
$$\phi_n\co \Gamma_B \rightarrow \ZZ_n,$$
where here we are using the identifications $H_1(B_0) \equiv \Gamma_B
\equiv H_1(B_1)$.  There are arbitrarily large pairs of such basis
elements; thus we may choose them so that the result of hyperbolic
Dehn filling $B_0$ and $B_1$ using these gives a two cusped hyperbolic
manifold with metric on the thick part as close to that of $X$ as
desired.

Now perform a very large Dehn filling (thus not distorting the
geometry of the thick part appreciably) along $A_1$.  We claim we
obtain $M$ with all the required properties.

For suppose that the geodesic representative for the unknotting tunnel
of $M$ had length at most $R$.  Let $T$ be the torus that forms the
boundary of a maximal horoball neighborhood of the cusp of $M$.  The
basepoint $m$ of $M$ lies on $T$.  We may pick $R$ large enough so
that $\pi_1(T,m)$ is generated by two curves of length at most $R$.
In addition, the geodesic representative for the unknotting tunnel can
be closed up to form a loop based at $m$ with length at most $2R$.
These three loops generate $\pi_1(M,m)$.  By construction, $H_1(M)
\cong {\mathbb Z} \oplus {\mathbb Z}_n$.  The image of $H_1(T)$ in
$H_1(M)$ is the first summand, and the image of the third loop is a
proper subgroup of the second summand.  Thus, these three loops cannot
generate $H_1(M)$, which is a contradiction.  Hence, the geodesic
representative for the unknotting tunnel of $M$ has length more than
$R$.  Since $R$ was arbitrarily large, this establishes Theorem
\ref{thm:long}.

\bibliographystyle{hamsplain}
\bibliography{biblio.bib}

\end{document}